\documentclass[journal]{IEEEtran}

\ifx\pdfoutput\undefined
\usepackage{graphicx} \else
\usepackage[pdftex]{graphicx} \fi
\usepackage{subfigure}
\usepackage{amsmath}
\usepackage{dsfont}
\usepackage{amssymb,amsthm}
\usepackage{algorithm}
\usepackage{subfigure}
\usepackage{cuted}
\usepackage{epstopdf}
\usepackage{cite}
\usepackage{soul,color,bm}
\usepackage{setspace}
\usepackage{amsbsy}

\newtheorem{thm}{Theorem}

\usepackage{setspace}
\usepackage{lipsum} 
\usepackage{multicol}
\usepackage{multirow}
\usepackage{float}
\usepackage{graphicx}

\usepackage{etoolbox}
\makeatletter
\patchcmd{\@begintheorem}{\textit}{\textbf}{}{}
\makeatother
\ifCLASSINFOpdf
\else
\fi

\begin{document}

\title{\huge{On Correcting Errors in Existing Mathematical Approaches for UAV Trajectory Design Considering No-Fly-Zones}}

\author{Kanghyun Heo, Gitae Park, and Kisong Lee,~\IEEEmembership{Member,~IEEE}
\vspace{-3mm}

\thanks{The authors are with the Department of Information and Communication Engineering, Dongguk University, Seoul 04620, South Korea (e-mail: kslee851105@gmail.com).}

}


\maketitle

\begin{abstract}
Motivated by the fact that current mathematical methods for the trajectory design of an unmanned aerial vehicle (UAV) considering no-fly-zones (NFZs) cannot perfectly avoid NFZs throughout the entire continuous trajectory, this study introduces a new constraint that ensures the complete avoidance of NFZs. Moreover, we provide mathematical proof demonstrating that a UAV operating within the proposed constraints will never violate NFZs. Under the proposed constraint on NFZs, we aim to optimize the scheduling, transmit power, length of the time slot, and the trajectory of the UAV to maximize the minimum throughput among ground nodes without violating NFZs. To find the optimal UAV strategy from the non-convex optimization problem formulated here, we use various optimization techniques, in this case quadratic transform, successive convex approximation, and the block coordinate descent algorithm. Simulation results confirm that the proposed constraint prevents NFZs from being violated over the entire trajectory in any scenario. Furthermore, the proposed scheme shows significantly higher throughput than the baseline scheme using the traditional NFZ constraint by achieving a zero outage probability due to NFZ violations.
\end{abstract}

\begin{IEEEkeywords}
Unmanned aerial vehicle, communications, no-fly-zone, trajectory design, convex optimization
\end{IEEEkeywords}

\IEEEpeerreviewmaketitle
\vspace{-3mm}
\section{Introduction}

In recent times, remarkable progress in the design and manufacturing of affordable unmanned aerial vehicles (UAVs) has resulted in a surge in their application across various domains, such as military operations, surveillance missions, cargo delivery, and even for exploring inaccessible or concealed regions. Particularly noteworthy is the elevated utilization of UAVs in wireless communications systems, primarily due to certain exceptional features of UAVs, such as their high mobility, rapid deployment, and on-demand connectivity \cite{Hayat2016, Zeng2016}. In particular, the communication links connecting UAVs and ground nodes (GNs) usually operate in a line-of-sight (LoS) manner, which can increase the communication capacity due to the high channel gain \cite{Lin2018}. Accordingly, UAV strategies related to their transmit power and trajectory have been optimized to maximize the throughput \cite{Zeng16} or to minimize the consumed energy \cite{Zeng19} in UAV-enabled communication systems. 

Considering practical geometrical constraints, i.e., no-fly-zones (NFZs), that impose restrictions on UAV flights over specific areas such as military bases and civil aviation airports, several recent studies have explored trajectory planning while taking into account the necessity of avoiding these NFZs \cite{Li18,Li20,Gao19,Wen23}. In particular, resource allocation and trajectory design were jointly optimized to maximize the communication throughput \cite{Li18} and to provide secure communications \cite{Li20,Gao19}, respectively, in the presence of the cylindrical NFZs. Moreover, a joint trajectory and pick-up design for UAV-assisted parcel delivery under NFZ constraints was proposed in \cite{Wen23}. Previous research \cite{Li18,Li20,Gao19,Wen23} has solved the problem by mathematically modeling the constraints of cylindrical NFZs, but we found that this method does not allow for perfect avoidance of NFZs over the entire continuous trajectory. Motivated by these observations, we introduce a new constraint that ensures the avoidance of NFZs with complete mathematical proof. The contributions of our study can be summarized as follows.

1) We reveal the shortcomings of existing mathematical methods used for UAV trajectory design considering cylindrical NFZs. While these approaches can avoid NFZ violations at specific discrete positions of the UAV, they cannot guarantee that such violations will not occur over the entire continuous trajectory. To ensure that UAVs do not consecutively violate NFZs, we propose a new constraint that realizes the avoidance of NFZs by introducing the concept of expanded NFZs and mathematically prove that a UAV operating within the proposed constraint will never violate NFZs.

2) We build the problem that optimizes the scheduling, transmit power, length of the time slot, and trajectory of the UAV to maximize the minimum throughput among GNs without violating NFZs. To deal with the non-convexity of the formulated optimization problem, we apply the quadratic transform (QT) and successive convex approximation (SCA) to approximate the problem as convex for each optimization variable and use a block coordinate descent algorithm to solve the problem. 

3) Through simulations under various scenarios, we confirm that under the proposed constraint, a UAV does not violate  NFZ in any case, unlike the situation under the existing constraint. As a result, the proposed scheme can achieve a zero outage probability due to NFZ violations in any environment and can show higher throughput than the conventional scheme.


     


\vspace{-3mm}
\section{System Model and Problem Statement}

\begin{figure}[t]
\centering
\includegraphics[width=0.7\linewidth]{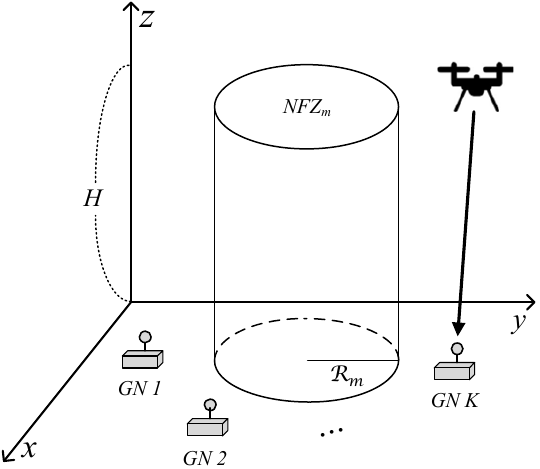} \vspace{-3mm} \caption{A UAV-enabled communication system.}\vspace{-3mm}
\label{fig1}
\end{figure}

Fig. \ref{fig1} shows a UAV-enabled communication system in which a UAV transmits data signals to a range of $K$ GNs. Let $T$ and $H$ denote the flight period and fixed altitude of the UAV, respectively, and $T$ is divided into $N$ time slots with different lengths $\delta[n]$ for $n \!\in\! \mathcal{N}\!=\!\{1,2,\cdots,N\}$. The constraints for $\delta[n]$ are then given by 
\begin{align}
    \delta[n] &> 0, ~~\forall n, \label{const1}\\
    \sum_{n\in\mathcal{N}}\delta[n]&= T. \label{const2}
\end{align}

The horizontal coordinates of the UAV at time slot $n$ are expressed as $\mathbf{q}[n]\!=\!(x[n],y[n])$, while the fixed horizontal coordinates of GN $k\!\in\! \mathcal{K}\!=\!\{1,2,\cdots,K\}$ are represented by $\mathbf{w}_k\!=\!(x_k,y_k)$. Moreover, the maximum flying speed of the UAV is denoted by $V$, which leads to a maximum flying distance of $V \delta[n]$ during time slot $n$. The UAV also needs to return to its initial location after one period to periodically support the GNs. Therefore, the constraints on UAV mobility are formulated as  
\begin{align}
\|\mathbf{q}[n]-\mathbf{q}[n-1]\| &\leq V \delta[n], ~~\forall n, \label{const3} \\
\mathbf{q}[0] &= \mathbf{q}[N]. \label{const4}
\end{align}

\begin{figure}[t]
  \begin{center}
    \subfigure[Conventional.]{
      \includegraphics[width=0.42\linewidth]{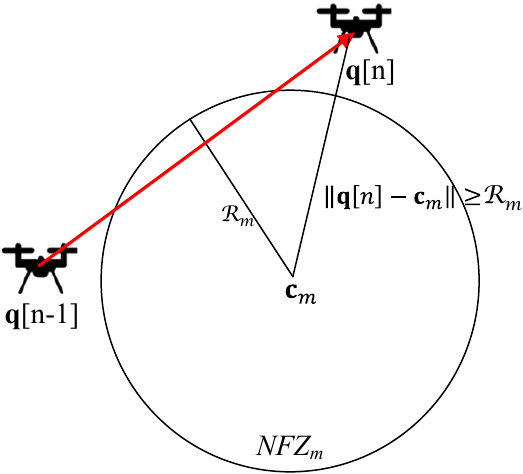}\label{fig2a}
    }
    \subfigure[Proposed.]{
      \includegraphics[width=0.48\linewidth]{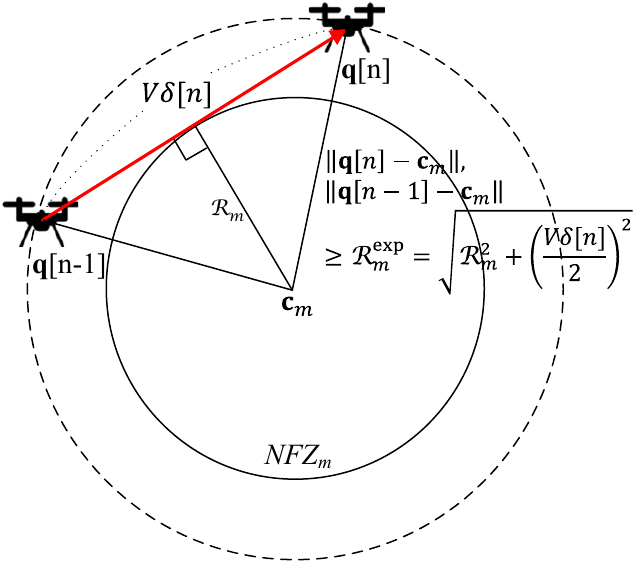}\label{fig2b}
    }
  \end{center} \vspace{-3mm}
\caption{Constraints on NFZ avoidance.}\vspace{-3mm}
\label{fig2}
\end{figure}

During its flight, the UAV may encounter $M$ non-overlapping NFZs with a cylindrical shape, which it must avoid. In earlier work \cite{Li18,Li20,Gao19,Wen23}, the constraint on the avoidance of NFZs is formulated as 
\begin{align}
\|\mathbf{q}[n]-\mathbf{c}_m\| &\geq \mathcal{R}_m, ~~ \forall m,~n, \label{const5}
\end{align}
where $\mathbf{c}_m=(x_m,y_m)$ and $\mathcal{R}_m$ are correspondingly the coordinate center and radius of the NFZ $m \!\in\! \mathcal{M}\!=\!\{1,2,\cdots,M\}$. However, as shown in Fig. \ref{fig2a}, the constraint \eqref{const5} only ensures that the UAV does not violate the NFZ at certain discrete points, i.e., $\mathbf{q}[n-1]$ and $\mathbf{q}[n]$, but it does not fully guarantee that the UAV will not violate the NFZ throughout the continuous trajectory between these two points. 

Inspired by this limitation of the previous study, we propose the following theorem to construct a new constraint that does not violate the NFZ throughout the continuous trajectory. 

\begin{thm} If the length of a line segment connecting any two points that are not interior points of a circle of radius $\mathrm{R}_B$ is less than $L$ and $L \leq 2\mathrm{R}_B$ holds, the distance between the line segment and the center of the circle is greater than or equal to $\sqrt{\mathrm{R}_B^2 - \left(\frac{L}{2}\right)^2}$. 
\end{thm}

\textit{Proof}: Please refer to the Appendix. \qed
\vspace{1mm}

In UAV communications as considered here, the trajectory of the UAV must not violate NFZs for all consecutive times. This means that the distance between the position of the UAV and the center of NFZ $m$, $\mathbf{c}_{m}$, must be greater than or equal to the radius of the NFZ $m$, $\mathcal{R}_m$, not only at $\mathbf{q}[n]$, $\forall n$ but also when transitioning from $\mathbf{q}[n-1]$ to $\mathbf{q}[n]$. To ensure this, we consider an expanded NFZ $m$ with radius $\mathcal{R}^{\textrm{exp}}_m[n]$ at time slot $n$. Given that the line segment $\overline{\mathbf{q}[n]\mathbf{q}[n-1]}$ is less than or equal to $V\delta[n]$, the distance between $\overline{\mathbf{q}[n]\mathbf{q}[n-1]}$ and $\mathbf{c}_{m}$ has a lower bound of $\mathcal{R}_m^{\textrm{LB}}=\sqrt{\left(\mathcal{R}^{\textrm{exp}}_m[n]\right)^2 -\left(\frac{V\delta[n]}{2}\right)^2}$ according to \textit{Theorem 1}. By setting the radius of the original NFZ to this lower bound, i.e., $\mathcal{R}_m=\sqrt{\left(\mathcal{R}^{\textrm{exp}}_m[n]\right)^2 -\left(\frac{V\delta[n]}{2}\right)^2}$, the distance between any consecutive point on the UAV trajectory and the center of the NFZ cannot be less than $\mathcal{R}_m$, which means that the UAV can completely avoid the NFZ for the entire consecutive trajectory. Therefore, we can determine the radius of the expanded NFZ $m$ at time slot $n$ as follows: 
\begin{align}
\mathcal{R}^{\textrm{exp}}_m[n] = \sqrt{\mathcal{R}_m^2+\left(\frac{V\delta[n]}{2}\right)^2}, ~~ \forall n. \label{radius1}
\end{align}

As shown in Fig. \ref{fig2b}, two consecutive time slots are directly related to the length of time slot; i.e., $\delta[n]$ is associated with $\mathbf{q}[n-1]$ and $\mathbf{q}[n]$, meaning that the distance between $\mathbf{c}_{m}$ and the UAV at points $\mathbf{q}[n-1]$ and $\mathbf{q}[n]$ must be greater than the radius of the expanded NFZ $m$, $\mathcal{R}^{\textrm{exp}}_m[n]$, to avoid a violation of the NFZ. Accordingly, we can modify the constraint \eqref{const5} using \eqref{radius1}, as follows:
\begin{align}
\|\mathbf{q}[n]-\mathbf{c}_{m}\| &\geq \sqrt{\mathcal{R}_m^2+\left(\frac{V\delta[n]}{2}\right)^2}, \nonumber\\
\|\mathbf{q}[n-1]-\mathbf{c}_{m}\| &\geq \sqrt{\mathcal{R}_m^2+\left(\frac{V\delta[n]}{2}\right)^2}, ~~ \forall m,~n. \label{const6} 
\end{align}
It should be noted that the UAV does not violate NFZs even when transitioning from $\mathbf{q}[n-1]$ to $\mathbf{q}[n]$ under the proposed constraint \eqref{const6}, as shown in Fig. \ref{fig2b}.

With $p[n]$ representing the transmit power of the UAV at time slot $n$, the UAV then has the following power constraints.
\begin{align}
\frac{1}{T}\sum_{n\in \mathcal{N}}\delta[n]p[n] &\leq P_{\textrm{avg}},  \label{const7}\\
0 \leq p[n] &\leq P_{\textrm{peak}}, ~~\forall n, \label{const8}
\end{align}
where $P_{\textrm{avg}}$ and $P_{\textrm{peak}}$ are the average and peak power budgets for the UAV, respectively.

We define a binary variable $s_k[n]$ to represent scheduling. If GN $k$ is served by the UAV at time slot $n$, then $s_k[n]=1$; otherwise $s_k[n]=0$. Furthermore, the UAV can only serve at most one GN at each time slot. Hence, the following constraints on $s_k[n]$ can be devised.
\begin{align}
    s_k[n] &\in \{0,1\}, ~~\forall k,~n,\label{const9} \\
    \sum_{k\in \mathcal{K}} s_k[n] &\leq 1, ~~\forall n. \label{const10}
\end{align}

We use the free-space path-loss model to account for the dominance of LoS links on the air-to-ground wireless channels, assuming that the Doppler effect due to the mobility of the UAV is fully compensated for in the GNs \cite{Zeng16,Li18,Li20,Gao19}. In consequence, the channel gain between the UAV and GN $k$ at time slot $n$ is given by
\begin{align}\label{channel}
h_{k}[n]=\frac{\beta_0}{d_{k}^2[n]}=\frac{\beta_0}{\|\mathbf{q}[n]-\mathbf{w}_k\|^2+H^2}, ~~\forall k,~n,
\end{align}
where $\beta_0$ is the channel power gain at the reference distance of $1$ m and $d_{k}[n]$ is the distance between the UAV and GN $k$ at time slot $n$. 

The achievable rate from the UAV to GN $k$ at time slot $n$ is then formulated as
\begin{align}\label{SE1}
r_{k}[n]\!&=\!\delta[n]\log_{2}\!\left(\! 1 \!+\! \frac{h_{k}[n]p[n]}{\sigma^2} \!\right) \nonumber\\
\!&=\!\delta[n]\log_{2}\!\left(\! 1 \!+\! \frac{\beta_0p[n]}{\sigma^2(\|\mathbf{q}[n]\!-\!\mathbf{w}_k\|^2\!+\!H^2)} \!\right),
\end{align}
where $\sigma^2$ represents the power of the additive white Gaussian noise.  
The time-averaged rate from the UAV to GN $k$ is also given by
\begin{align}\label{SE2}
\bar{r}_k&=\frac{1}{T}\sum_{n\in \mathcal{N}}s_k[n]r_{k}[n], ~~\forall k.
\end{align}

In this study, our goal is to maximize the minimum average rate among the GNs while avoiding multiple NFZs throughout the continuous trajectory. 
To achieve this, we formulate the problem that optimizes the scheduling $\mathbf{S}\triangleq\{s_k[n], ~\forall k,~n\}$, transmit power $\mathbf{P}\triangleq\{p[n], ~\forall n\}$, trajectory $\mathbf{Q}\triangleq\{\mathbf{q}[n], ~\forall n\}$, and length of the time slot $\mathbf{\Delta}\triangleq\{\delta[n], ~\forall n\}$, as follows: 
\begin{align} 
\textbf{(P0):} ~\max_{\mathbf{S},~{\mathbf{P}},~\mathbf{Q},~\mathbf{\Delta}} & ~~ \min_{k\in \mathcal{K}} \bar{r}_k \nonumber  \\
\textrm{s. t.} ~~~ & ~~~\eqref{const1}-\eqref{const4},~\eqref{const6}-\eqref{const10}. \nonumber
\end{align}
The optimization problem $\textbf{(P0)}$ is a non-convex mixed-integer program because $\mathbf{S}$ is a binary variable and the objective function is not jointly concave with respect to (w.r.t.) the optimization variables, e.g., $\mathbf{S}$, $\mathbf{P}$, $\mathbf{Q}$ and $\mathbf{\Delta}$. Therefore, it is difficult analytically to derive a globally optimal solution.

\vspace{-3mm}
\section{Proposed Algorithm}

Problem $\textbf{(P0)}$ is challenging to optimize due to its non-convexity, and for this reason, we divide the original problem into three subproblems. We then use the QT and SCA to make each subproblem convex for each optimization variable and solve the problem using an existing convex solver, in this case, CVX \cite{Grant}.

\vspace{-3mm}
\subsection{Scheduling}

Introducing an auxiliary variable $r_{\textrm{min}}$, the problem to find the optimal $\mathbf{S}$ for fixed values of $\mathbf{P}$, $\mathbf{Q}$, and $\mathbf{\Delta}$ can be formulated as follows:
\begin{align} 
\textbf{(P1):} ~\max_{\mathbf{S},~r_{\textrm{min}}} & ~~ r_{\textrm{min}} \nonumber  \\
\textrm{s. t.} & ~~~\bar{r}_k \geq r_{\textrm{min}}, ~~\forall k, \label{const15} \\
& ~~~\eqref{const9},~\eqref{const10}. \nonumber
\end{align}

To make problem $\textbf{(P1)}$ more tractable, first we relax the binary variable $s_k[n]$ in \eqref{const9} into a continuous variable. To preserve the binary nature of $s_k[n]$, we also consider the following constraints.
\begin{align}
    &s_k[n] \in [0,1],~~\forall k,~n, \label{const9-1}\\
    &s_k[n](1-s_k[n]) \leq 0,~~\forall k,~n. \label{sch} 
\end{align}
Furthermore, to make the constraint \eqref{sch} a convex set, we apply the first-order Taylor expansion to find the upper bound of $s_k[n](1-s_k[n])$ and use the penalty convex-concave procedure (CCP) by introducing a slack variable $\phi_k[n] \geq 0$ \cite{Lipp16}, as follows. 
\begin{align}
    s_k[n](1-2s_k^r[n])+(s_k^{r}[n])^2 \leq \phi_k[n],~~\forall k,~n, \label{const10-1} 
\end{align}
where $s_k^{r}[n]$ is a scheduling indicator for the $r$-th iteration and $\phi_k[n]$ acts as a penalty to enlarge the initial feasible region of $s_k[n]$ for stable optimization. 

Using $\eqref{const9-1}$ and $\eqref{const10-1}$, problem $\textbf{(P1)}$ can be transformed into the following convex optimization problem.
\begin{align} 
\textbf{(P1-2):} ~\max_{\mathbf{S},~\mathbf{\Phi} \succeq 0,~r_{\textrm{min}}} & ~~ r_{\textrm{min}} -  \eta \sum_{k\in \mathcal{K}}\sum_{n\in \mathcal{N}}\phi_k[n] \nonumber  \\
\textrm{s. t.}~~~ & ~~~\eqref{const10},~\eqref{const15},~\eqref{const9-1},~\eqref{const10-1}, \nonumber
\end{align}
where $\mathbf{\Phi}\triangleq\{\phi_k[n], ~\forall k,~n\}$ and $\eta > 0$ is a regularization factor that adjusts the influence of the penalty term $\sum_{k\in \mathcal{K}}\sum_{n\in \mathcal{N}}\phi_k[n]$, which controls the feasibility of the constraint. It should be noted that problem $\textbf{(P1-2)}$ is optimized to maximize $r_{\textrm{min}}$ at low $\eta$ values, and as $\eta$ increases, it is optimized to maximize $r_{\textrm{min}}$ while remaining consistent with the binary nature of $\mathbf{S}$ \cite{Lipp16}. 

\vspace{-3mm}
\subsection{Power Allocation}

For fixed values of $\mathbf{S}$, $\mathbf{Q}$, and $\mathbf{\Delta}$, the problem to find the optimal $\mathbf{P}$ can be formulated as follows:
\begin{align} 
\textbf{(P2):} ~\max_{\mathbf{P},~r_{\textrm{min}}} & ~~ r_{\textrm{min}} \nonumber  \\
\textrm{s. t.} & ~~~\eqref{const7},~\eqref{const8},~\eqref{const15}. \nonumber
\end{align}
It is obvious that problem $\textbf{(P2)}$ is a convex problem that can easily be solved by CVX. 

\vspace{-3mm}
\subsection{Trajectory and Length of the Time Slot}

For fixed values of $\mathbf{S}$ and $\mathbf{P}$, the problem of jointly finding the optimal $\mathbf{Q}$ and $\mathbf{\Delta}$ can be formulated as follows:
\begin{align} 
\textbf{(P3):} ~\max_{\mathbf{Q},~\mathbf{\Delta},~r_{\textrm{min}}} & ~~ r_{\textrm{min}} \nonumber  \\
\textrm{s. t.} ~~ & ~~~\eqref{const1}-\eqref{const4},~\eqref{const6},~\eqref{const15}. \nonumber
\end{align}
In problem $\textbf{(P3)}$, the constraints $\eqref{const1}$-$\eqref{const4}$ are convex sets but $\eqref{const6}$ and $\eqref{const15}$ are not. 

To make $\eqref{const6}$ a convex set, we apply the first-order Taylor expansion to find the lower bound of $\|\mathbf{q}[n]\!-\!\mathbf{c}_{m}\|^2$, as follows.
\begin{align}
    \|\mathbf{q}[n]\!-\!\mathbf{c}_{m}\|^2 &\geq
    2(\mathbf{q}^r[n]\!-\!\mathbf{c}_m)^T(\mathbf{q}[n]\!-\!\mathbf{q}^r[n])\!+\!\|\mathbf{q}^r[n]\!-\!\mathbf{c}_{m}\|^2  \nonumber \\ 
    &\triangleq G^{\textrm{LB}}_m[n] \label{glb}
\end{align}
Using $\eqref{glb}$, the constraint $\eqref{const6}$ can be transformed into the following convex set. 
\begin{align}
    G^{\textrm{LB}}_m[n] &\geq \mathcal{R}_m^2\!+\!\left(\frac{V\delta[n]}{2}\right)^2, \nonumber\\
    G^{\textrm{LB}}_m[n\!-\!1] &\geq \mathcal{R}_m^2\!+\!\left(\frac{V\delta[n]}{2}\right)^2, ~~ \forall m,~n. 
\label{const6-1}
\end{align}

To deal with the non-convexity of $r_k[n]$ in $\eqref{const15}$, first we transform $r_k[n]$ into its equivalent form. 
\begin{align}
    r_k[n] \!&=\! \delta[n]\log_2\!\left(\sigma^2\!\left(\|\mathbf{q}[n] \!-\! \mathbf{w}_k\|^2 \!+\! H^2\right) \!+\! \beta_0p[n]\right) \nonumber \\ 
    &~~~~ -\! \delta[n]\log_2\!\left(\sigma^2\!\left(\|\mathbf{q}[n] \!-\! \mathbf{w}_k\|^2 \!+\! H^2\right)\!\right) \nonumber \\ 
    \!&=\!  \left(\!\frac{\log_2\!\left(\sigma^2\!\left(\|\mathbf{q}[n] \!-\! \mathbf{w}_k\|^2 \!+\! H^2\right) \!+\! \beta_0p[n]\right)}{1/\delta[n]}\!\!\right) \nonumber \\
    &~~~~+\!\! \left(\!\frac{\log_2\!\left(\!\frac{1}{\sigma^2(\|\mathbf{q}[n] - \mathbf{w}_k\|^2 + H^2)}\!\right)}{1/\delta[n]}\!\!\right). \label{rkn}
\end{align}

To make the logarithm terms in \eqref{rkn} concave w.r.t. $\mathbf{q}[n]$, we use the first-order Taylor expansion, as follows: 
\begin{align}
    &\|\mathbf{q}[n] \!-\! \mathbf{w}_k\|^2 \!+\! H^2  \geq 2(\mathbf{q}^r[n] \!-\! \mathbf{w}_k)^T(\mathbf{q}[n] \!-\! \mathbf{q}^r[n]) \nonumber\\ 
    &~~~~~~~~~~~~~~~~~~~~~~~~+\|\mathbf{q}^r[n] \!-\! \mathbf{w}_k\|^2 \!+\! H^2 \triangleq A^{\textrm{LB}}_k[n],
    \label{const_A} 
\end{align}
\begin{align}
    &\frac{1}{\|\mathbf{q}[n] \!-\! \mathbf{w}_k\|^2 \!+\! H^2} \geq -\frac{\|\mathbf{q}[n] \!-\! \mathbf{w}_k\|^2 \!-\! \|\mathbf{q}^r[n] \!-\! \mathbf{w}_k\|^2}{(\|\mathbf{q}^r[n] \!-\! \mathbf{w}_k\|^2 \!+\! H^2)^2} \nonumber\\ 
    & ~~~~~~~~~~~~~~~~~~~~~~~~+\frac{1}{\|\mathbf{q}^r[n] \!-\! \mathbf{w}_k\|^2 \!+\! H^2} \triangleq B^{\textrm{LB}}_k[n].
    \label{const_B}
\end{align}
Using $\eqref{const_A}$ and $\eqref{const_B}$, the lower bound of $r_k[n]$ can be formulated as 
\begin{align}
    r_k[n] &\geq\!  \left(\!\frac{\log_2\!\left(\sigma^2\!A^{\textrm{LB}}_k[n] \!+\! \beta_0p[n]\right)\!+\!\gamma_1}{1/\delta[n]}\!\!\right) \!\!+\!\! \left(\!\frac{\log_2\!\left(\!\frac{B^{\textrm{LB}}_k[n]}{\sigma^2}\!\right)\!+\!\gamma_2}{1/\delta[n]}\!\!\right)\nonumber\\
    &~~~~~-\delta[n](\gamma_1 \!+\! \gamma_2) \triangleq r^{\textrm{LB}}_k[n], \label{rkn2}
\end{align}
where $\gamma_1$ and $\gamma_2$ are constants required to make the QT possible, which will be explained in \eqref{tranformed}.
Moreover, the lower bounds $A^{\textrm{LB}}_k[n]$ and $B^{\textrm{LB}}_k[n]$ must satisfy the following constraints for the feasibility of \eqref{rkn2}. 
\begin{align}
    A^{\textrm{LB}}_k[n]&> \frac{-\beta_0p[n]}{\sigma^2}, ~~ \forall k,~n,  \label{A1}\\
    B^{\textrm{LB}}_k[n]&> 0, ~~ \forall k,~n. \label{B1}
\end{align}

Because the first and second terms in \eqref{rkn2} have a concave-convex form w.r.t. $\mathbf{q}[n]$ and $\delta[n]$, the QT can be used to transform $r^{\textrm{LB}}_k[n]$ into the following equivalent form \cite{Shen18}. 
\begin{align}
    &f^{\textrm{LB}}_k[n] \triangleq  2\lambda_k[n]\!\sqrt{\log_2\!\left(\sigma^2A^{\textrm{LB}}_k[n]\!+\!\beta_0p[n] \right) \!+\!\gamma_1}\!-\! \frac{\lambda_k^2[n]}{\delta[n]} \nonumber \\
    & +2\mu_k[n]\!\sqrt{\log_2\!\left(\!\frac{B^{\textrm{LB}}_k[n]}{\sigma^2}\!\right) \!+\! \gamma_2} \!-\! \frac{\mu_k^2[n]}{\delta[n]} \!-\!\delta[n](\gamma_1 \!+\! \gamma_2),
    \label{tranformed}
\end{align}
where $\gamma_1$ and $\gamma_2$ ensure the feasibility of $f^{\textrm{LB}}_k[n]$ by preventing the square root function from going negative, which is why they are added to \eqref{rkn2}. In addition, $\pmb{\lambda}=\{\lambda_k[n], ~\forall k,~n\}$ and $\pmb{\mu}=\{\mu_k[n], ~\forall k,~n\}$ denote the set of non-negative auxiliary variables for the QT. In \eqref{tranformed}, the terms related to $\mathbf{q}[n]$, e.g., $2\lambda_k[n]\!\sqrt{\log_2\!\left(\sigma^2A^{\textrm{LB}}_k[n]\!+\!\beta_0p[n] \right) \!+\!\gamma_1}$ and $2\mu_k[n]\!\sqrt{\log_2\!\left(\!\frac{B^{\textrm{LB}}_k[n]}{\sigma^2}\!\right) \!+\! \gamma_2}$, are concave w.r.t. $\mathbf{q}[n]$ because the square root is a concave and non-decreasing function and the logarithm is a concave function. Furthermore, the terms related to $\delta[n]$, e.g., $-\frac{\lambda_k^2[n]}{\delta[n]}$, $-\frac{\mu_k^2[n]}{\delta[n]}$ and $-\delta[n](\gamma_1 \!+\! \gamma_2)$, are concave w.r.t. the positive value of $\delta[n]$. Therefore, $f^{\textrm{LB}}_k[n]$ is jointly concave w.r.t. $\mathbf{q}[n]$ and $\delta[n]$ when the auxiliary variables $\lambda_k[n]$ and $\mu_k[n]$ are fixed. 

Furthermore, because $f^{\textrm{LB}}_k[n]$ is concave w.r.t. $\lambda_k[n]$ and $\mu_k[n]$ for fixed values of $\mathbf{q}[n]$ and $\delta[n]$, we can find the optimal values of $\lambda_k[n]$ and $\mu_k[n]$ from $\frac{\partial f^{\textrm{LB}}_k[n]}{\partial \lambda_k[n]}=0$ and $\frac{\partial f^{\textrm{LB}}_k[n]}{\partial \mu_k[n]}=0$, respectively, as follows:
\begin{align}
    \lambda_k^*[n] &= \delta[n]\sqrt{\log_2\!\left(\sigma^2A^{\textrm{LB}}_k[n] \!+\! \beta_0p[n]\right)\!+\!\gamma_1}, \label{lambda}\\
    \mu_k^*[n] &= \delta[n]\sqrt{\log_2\!\left(\!\frac{B^{\textrm{LB}}_k[n]}{\sigma^2}\!\right) \!+\! \gamma_2}. \label{mu}
\end{align}

The time-averaged value of $f^{\textrm{LB}}_k$ can then be represented by
\begin{align}
    \bar{f}^{\textrm{LB}}_k = \frac{1}{T}\sum_{n \in \mathcal{N}}s_k[n]f^{\textrm{LB}}_k[n]. \label{F}
\end{align}
Using \eqref{F}, the constraint $\eqref{const15}$ can be transformed into the following convex set: 
\begin{align}
    \bar{f}^{\textrm{LB}}_k \geq r_{\textrm{min}}, ~~ \forall k. \label{const15-1}
\end{align}

Using \eqref{const6-1} and \eqref{const15-1}, problem $\textbf{(P3)}$ can be reformulated as the following convex problem.
\begin{align} 
\textbf{(P3-2):} ~\max_{\mathbf{Q},~\mathbf{\Delta},~\pmb{\lambda},~\pmb{\mu},~r_{\textrm{min}}} & ~~~ r_{\textrm{min}} \nonumber  \\
\textrm{s. t.} ~~~~~ & ~~~\eqref{const1}-\eqref{const4},~\eqref{const6-1},~\eqref{A1},~\eqref{B1},~\eqref{const15-1}.\nonumber 
\end{align}

To handle the non-convexity of problem $\textbf{(P0)}$, we develop three subproblems, each convex for each optimization variable, and then iteratively solve these subproblems using a convex solver until convergence. Algorithm \ref{Alg1} summarizes the detailed procedure of the proposed method.

\begin{algorithm}[ht]
    \caption{Proposed Algorithm} \label{Alg1}\small 
    1:$~~$Set $r\!=\!1$ \\
    2:$~~$Initialize $\mathbf{S}^{r}$, $\mathbf{P}^{r}$, $\mathbf{Q}^{r}$, $\mathbf{\Delta}^{r}$ and $\eta^r$ \\
    3:$~~$\textbf{repeat}  \\
    4:$~~~~$Update $\{\pmb{\lambda},\pmb{\mu}\}$ by \eqref{lambda} and \eqref{mu}, respectively  \\
    5:$~~~~$Find $\{\mathbf{Q}^{r+1}\!,\mathbf{\Delta}^{r+1}\}$ by solving $\textbf{(P3-2)}$ for given $\{\mathbf{S}^{r}\!,\mathbf{P}^{r}\}$     \\
    6:$~~~~$Find $\mathbf{S}^{r+1}$ by solving $\textbf{(P1-2)}$ for given $\{\mathbf{P}^{r}\!,\mathbf{Q}^{r+1}\!,\mathbf{\Delta}^{r+1}\!\}$ \\
    7:$~~~~$Update $\eta^{r+1} = \min \{\kappa\eta^r, \eta_{\textrm{max}}\}$ \\
    8:$~~~~$Find $\mathbf{P}^{r+1}$ by solving $\textbf{(P2)}$ for given $\{\mathbf{S}^{r+1}\!,\mathbf{Q}^{r+1}\!,\mathbf{\Delta}^{r+1}\!\}$ \\
    9:$~~~~$Update $r \leftarrow r+1$  \\
    10:$~$\textbf{until} The increase in the objective value is less than $\epsilon > 0$.
\end{algorithm}

\begin{figure*}[t]
  \begin{center}
    \subfigure[Trajectory and scheduling.]{
      \includegraphics[width=0.45\linewidth]{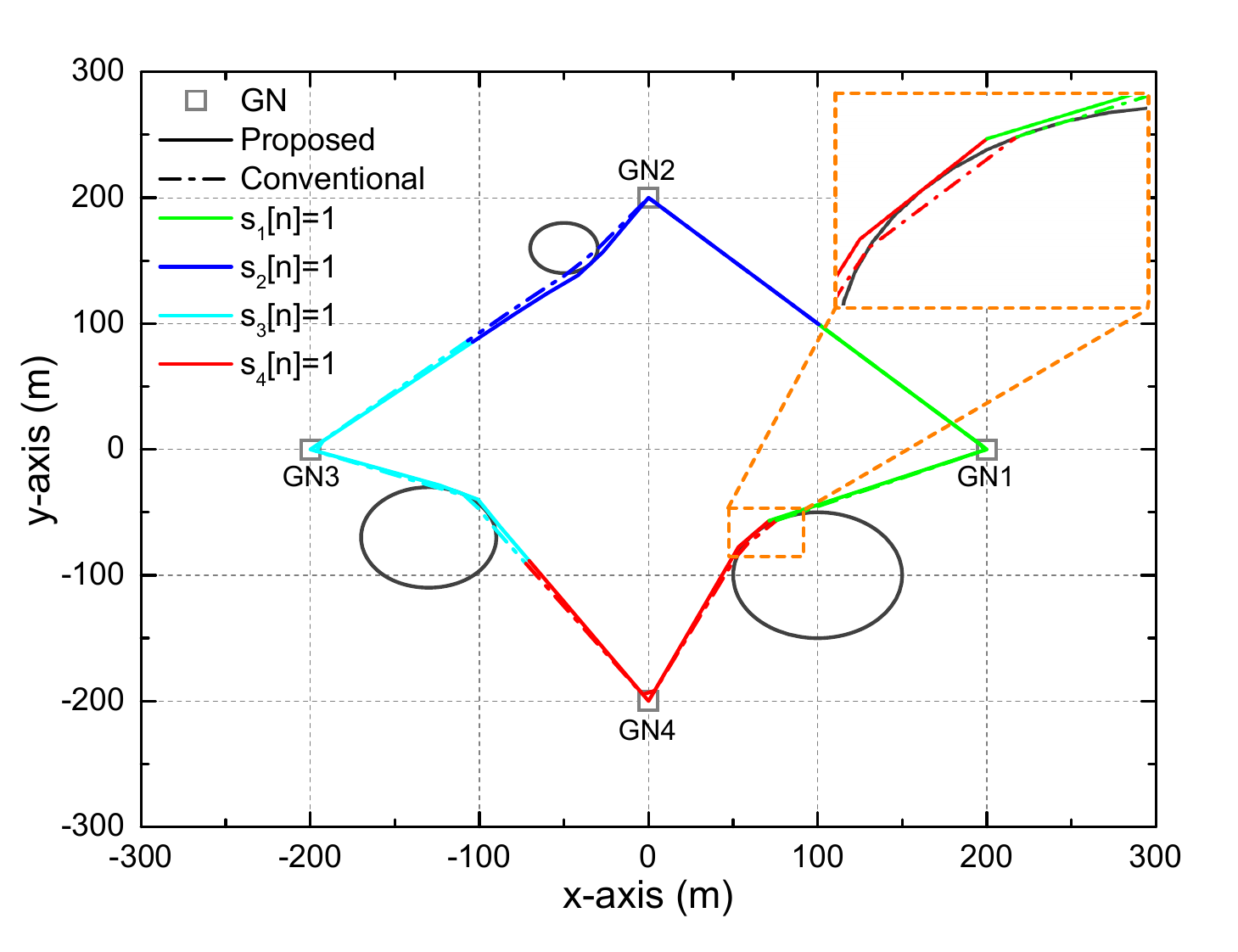}\label{R1-1}
    }
    \subfigure[Length of time slots and transmit power.]{
      \includegraphics[width=0.45\linewidth]{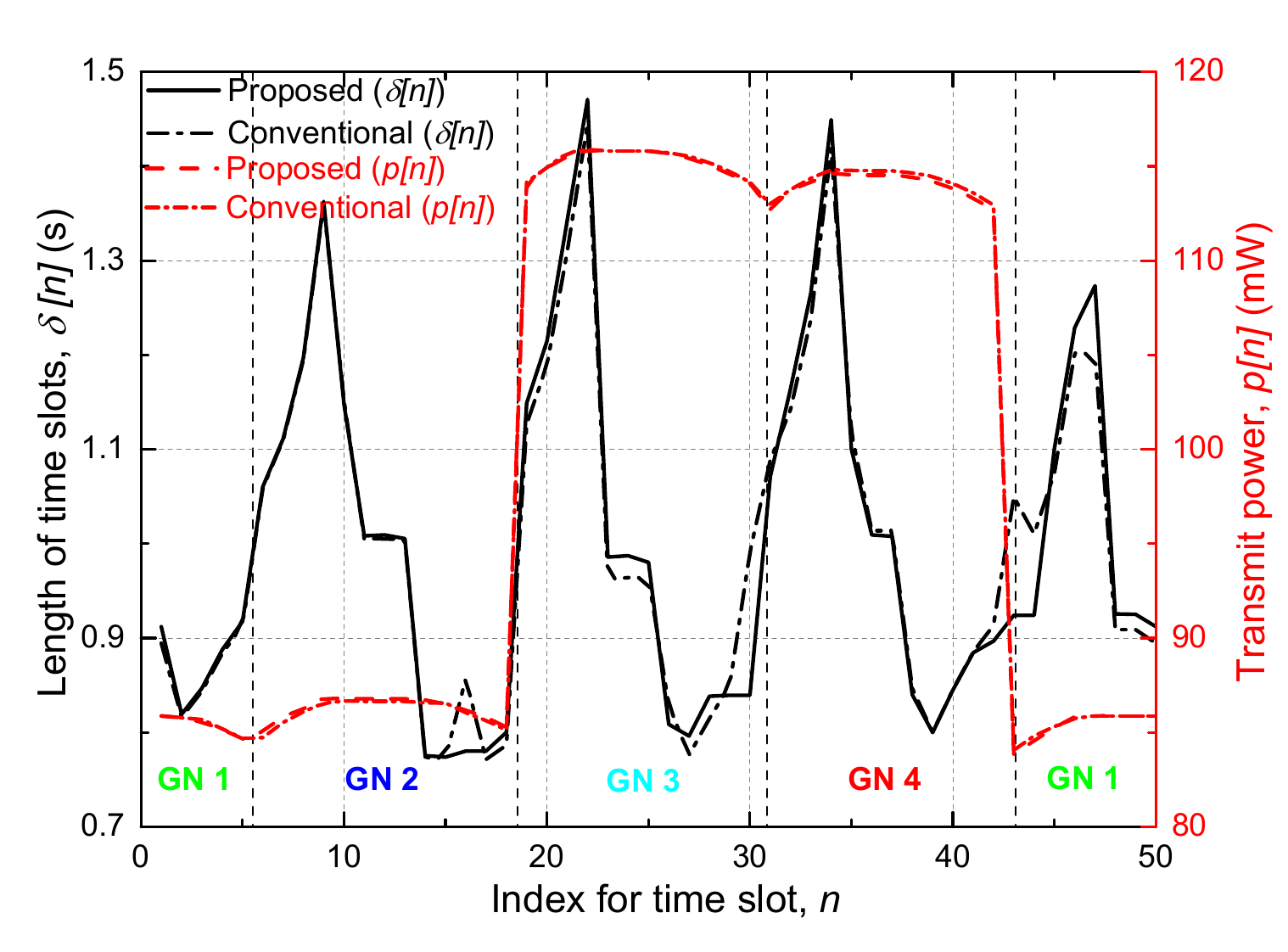}\label{R1-2}
    }
  \end{center} 
\caption{Comparison of resource allocation between proposed and conventional schemes.}
\label{R1}
\end{figure*}

\begin{figure*}[t]
  \begin{center}
    \subfigure[Average spectral efficiency.]{
      \includegraphics[width=0.45\linewidth]{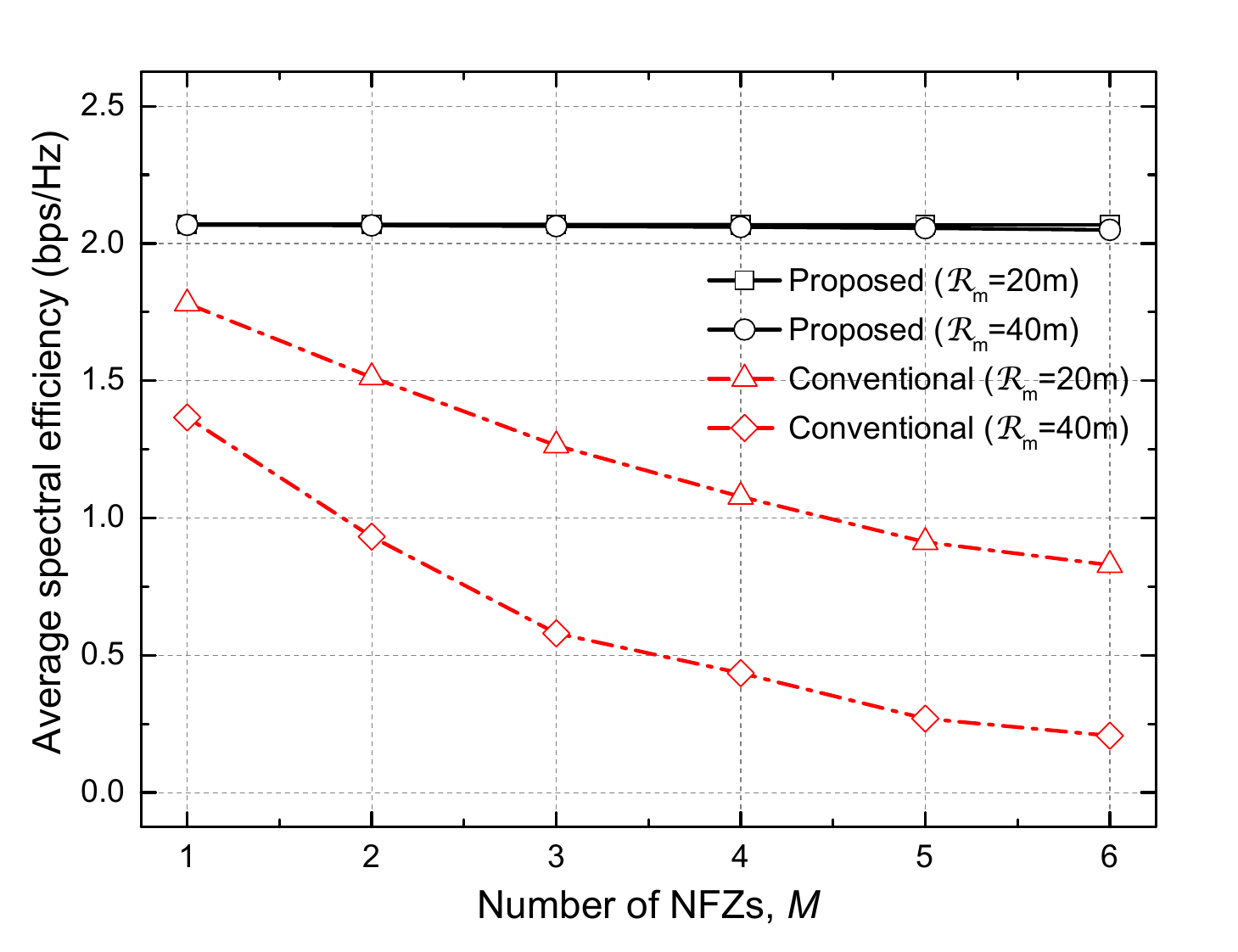}\label{R2-1}
    }
    \subfigure[Outage probability.]{
      \includegraphics[width=0.45\linewidth]{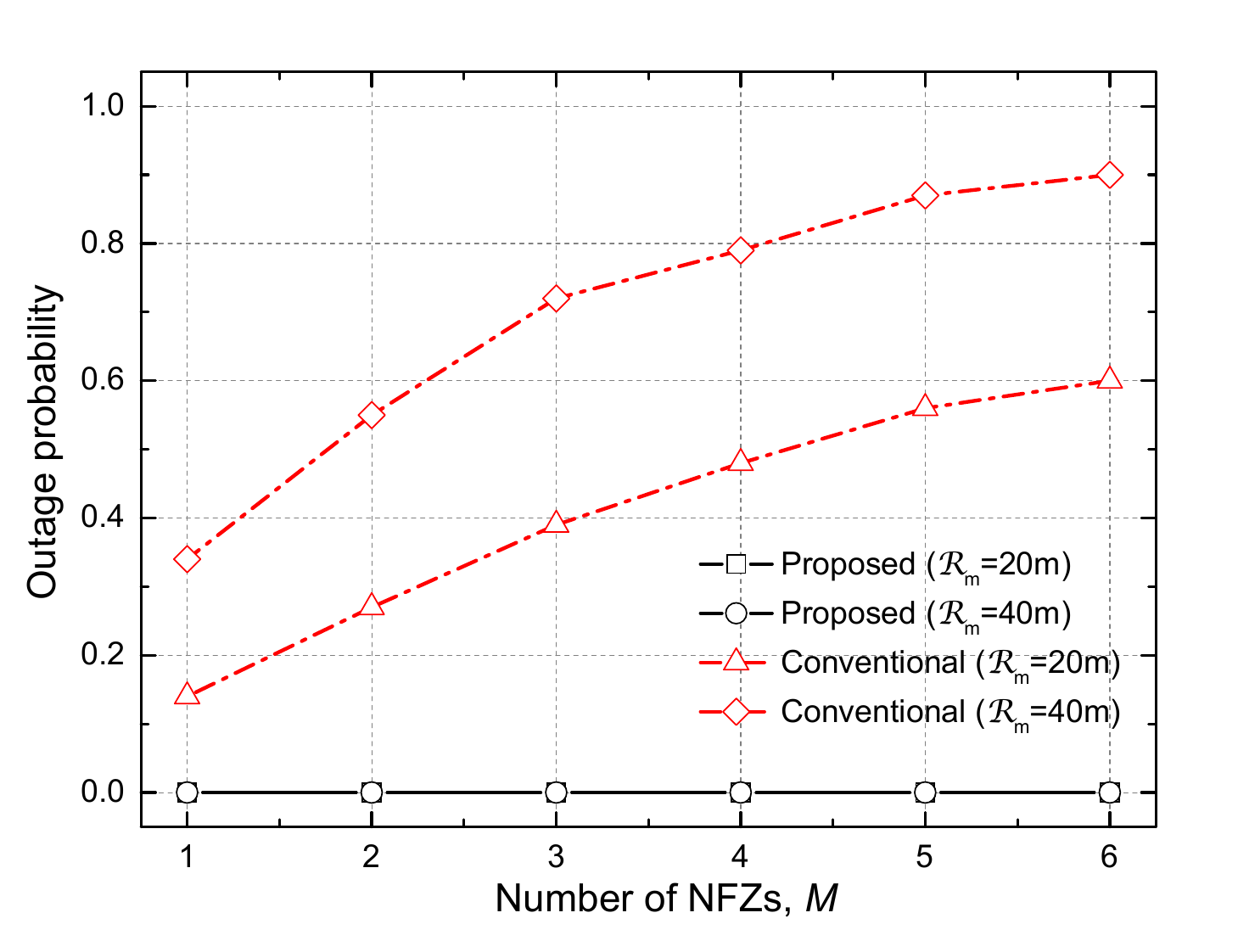}\label{R2-2}
    }
  \end{center} 
\caption{Performance comparison between the proposed and conventional schemes for the radius of the NFZ and the number of NFZs.}
\label{R2}
\end{figure*}

The convergence of Algorithm \ref{Alg1} can be guaranteed because the objective function of each subproblem is a non-decreasing function at each iteration and is bounded to a finite value. Moreover, the computational complexity of the proposed algorithm is $\emph{O}\!\left(R_C(KN)^{3.5}\log(1/\epsilon)\right)$, where $R_C$ is the number of iterations required for convergence (lines 3--10). Indeed, the proposed algorithm possesses the computational complexity of polynomial time, making it well suited for real-time implementations \cite{Li20,Leiserson}.

\vspace{-3mm}
\section{Results and Discussions}

For performance evaluations, we randomly generate GNs and NFZs in a square area $400$ m $\times$ $400$ m in size. The following parameters are considered as default values \cite{Zeng16,Zeng19,Li18,Li20,Gao19,Wen23}: $T\!=\!50$ s, $N\!=\!50$, $M\!=\!3$, $K\!=\!4$, $H\!=\!30$ m, $V_{\textrm{max}}\!=\!30$ m/s, $P_{\textrm{avg}}\!=\!20$ dBm, $P_{\textrm{peak}}\!=\!4P_{\textrm{avg}}$, $\beta_0\!=\!-30$ dB, $\sigma^2\!=\!-70$ dBm, $\gamma_1\!=\!\gamma_2\!=\!100$, $\eta^1\!=\!1$, $\eta_{\textrm{max}}\!=\!10^6$, $\kappa\!=\!1.2$ and $\epsilon \!=\! 10^{-4}$ . 

Fig. \ref{R1} shows a resource allocation comparison between the proposed scheme and the conventional scheme using the NFZ constraint \eqref{const5}. As shown in Fig. \ref{R1-1}, the UAV transmits a signal to the closest GN by setting the scheduling indicator to $1$ as it approaches each GN. Moreover, if the UAV encounters NFZs while flying, it will take a trajectory to avoid them. However, the UAV using the conventional scheme encroaches slightly into NFZs, as shown in the enlarged figure, while the UAV using the proposed scheme does not violate any NFZs. Fig. \ref{R1-2} shows the optimized lengths of the time slots ($\delta[n]$) and transmit power ($p[n]$) for both schemes. The UAV schedules GN 1, GN 2, GN 3, GN 4, and then GN 1 again. In this case, GN 1 and GN 2 are each allocated 13 time slots, and GN 3 and GN 4 are each allocated 12 time slots. Because GN 3 and GN 4 have one less time slot allocated for scheduling, the UAV allocates higher power when servicing GN 3 and GN 4 to compensate for this. We can also observe that the UAV increases $\delta[n]$ when it hovers directly over each GN to stay longer and serve more efficiently; e.g., $n=9,~22,~33,$ and $47$. In the proposed scheme, the UAV decreases $\delta[n]$ near the NFZ to get as close to the NFZ as possible to avoid it in a time-efficient manner, whereas in the conventional scheme, it increases $\delta[n]$ to pass through the NFZ quickly at once. This occurs because the traditional constraint \eqref{const5} focuses on avoiding NFZs at specific discrete points, not along the entire continuous trajectory. As a result, it violates NFZs near $n=16,~30,$ and $43$.

Fig. \ref{R2} shows a performance comparison between the proposed and conventional schemes with reference to the radius of the NFZ ($\mathcal{R}_m$) and the number of NFZs ($M$): (a) the average spectral efficiency and (b) the outage probability. Here, the outage probability is the probability that the UAV violates NFZs and the average spectral efficiency is set to zero when an outage occurs. As $\mathcal{R}_m$ and $M$ increase, the outage probability of the conventional scheme increases significantly because the UAV frequently violates NFZs. On the other hand, the proposed scheme shows stable performances that are not affected by the increases in $\mathcal{R}_m$ and $M$ because the UAV can perfectly avoid NFZs, as proven in \emph{Theorem 1}. This result demonstrates the effectiveness of the proposed method in terms of NFZ avoidance.

\vspace{-3mm}
\section{Conclusions}

In this paper, we proposed a new constraint supported by a rigorous mathematical proof that ensures the absolute avoidance of NFZs over an entire continuous trajectory of a UAV. In particular, we optimized the scheduling, transmit power, length of the time slot, and the trajectory of the UAV to maximize the minimum throughput among GNs without violating NFZs. To tackle the non-convexity of the optimization problem formulated, we applied the QT and SCA approaches to make the problem convex for each optimization variable and solved the problem using the block coordinate descent algorithm. The simulation results demonstrated that the UAV can efficiently provide services to GNs without violating NFZs under the proposed constraint while violating NFZs under the existing constraint. As a result, the proposed scheme can achieve higher throughput than the conventional scheme by significantly reducing the probability of outages due to NFZ violations. By rectifying the existing incorrect NFZ constraints, our study has the potential to be implemented in the designs of UAV trajectories in real-world scenarios that encompass multiple NFZs.

\vspace{-3mm}
\section*{Appendix}

Consider a circle $B$ with radius $\mathrm{R}_B$ and let $\mathcal{B}$ be the set of the interior points of $B$. We assume that the length of the line segment connecting any two points that are not contained in $\mathcal{B}$, i.e., $\mathbf{z}_1$, $\mathbf{z}_2$ $\notin \mathcal{B}$, is less than or equal to $L$, i.e., $\|\overline{\mathbf{z}_1\mathbf{z}_2}\| \leq L$. Then, there are two cases: when $\overline{\mathbf{z}_1\mathbf{z}_2}$ includes some point $\mathbf{b}$ contained in $\mathcal{B}$, i.e., $\overline{\mathbf{z}_1\mathbf{z}_2} \cap \mathcal{B} \neq \emptyset$ and when it does not, i.e., $\overline{\mathbf{z}_1\mathbf{z}_2} \cap \mathcal{B} = \emptyset$. 

When $\overline{\mathbf{z}_1\mathbf{z}_2}$ does not contain any $\mathbf{b} \in$ $\mathcal{B}$, it is obvious that the distance from any $\mathbf{z} \in \overline{\mathbf{z}_1\mathbf{z}_2}$ to the center of $B$, denoted by $\mathbf{c}_B$, is greater than or equal to $\mathrm{R}_B$. Therefore, the distance between the line segment $\overline{\mathbf{z}_1\mathbf{z}_2}$ and the center of the circle $\mathbf{c}_B$ is always greater than or equal to $\mathrm{R}_B$, as follows: 
\begin{align}
    d(\overline{\mathbf{z}_1\mathbf{z}_2},\mathbf{c}_B) = \min_{\mathbf{z}\in \overline{\mathbf{z}_1\mathbf{z}_2}} \|\mathbf{z}-\mathbf{c}_B \| \geq \mathrm{R}_B.
\end{align}


When $\overline{\mathbf{z}_1\mathbf{z}_2}$ contains any $\mathbf{b} \in$ $\mathcal{B}$, $\overline{\mathbf{z}_1 \mathbf{b}}$ and $\overline{\mathbf{b} \mathbf{z}_2}$ are also line segments which are the subsets of $\overline{\mathbf{z}_1\mathbf{z}_2}$, i.e., $ \overline{\mathbf{z}_1 \mathbf{b}}$, $\overline{\mathbf{b} \mathbf{z}_2}$ $\subset \overline{\mathbf{z}_1\mathbf{z}_2}$. The line segment from the interior point of the circle, e.g., $\mathbf{b}$, to the exterior point, e.g., $\mathbf{z}_1$ or $\mathbf{z}_2$, passes through the boundary of the circle. Let $\mathbf{z}_1'$ and $\mathbf{z}_2'$ denote the boundary points of the circle, i.e., $\mathbf{z}_1' \in \overline{\mathbf{z}_1\mathbf{b}}$, $\mathbf{z}_2' \in \overline{\mathbf{b}\mathbf{z}_2}$, respectively. Because $\overline{\mathbf{z}_1\mathbf{b}}~ \cap~ \overline{\mathbf{b}\mathbf{z}_2} = \{\mathbf{b}\}$, $\mathbf{z}_1' \neq \mathbf{z}_2'$ holds. Therefore, $\overline{\mathbf{z}_1'\mathbf{z}_2'}$ is the chord of the circle and is simultaneously a subset of $\overline{\mathbf{z}_1\mathbf{z}_2}$ such that $\overline{\mathbf{z}_1'\mathbf{z}_2'} \subset \overline{\mathbf{z}_1\mathbf{z}_2}$. Then, we can build the following relationship, $\|\overline{\mathbf{z}_1'\mathbf{z}_2'}\| \leq \|\overline{\mathbf{z}_1\mathbf{z}_2}\| \leq L$. 

Because any point on the chord is an interior point or boundary point of the circle, the following condition holds. 
\begin{align}
    \min_{\mathbf{z} \in \overline{\mathbf{z}_1'\mathbf{z}_2'}}\|\mathbf{z}-\mathbf{c}_B\| \leq \mathrm{R}_B < \min_{\mathbf{z} \in  \overline{\mathbf{z}_1\mathbf{z}_2}\backslash \overline{\mathbf{z}_1'\mathbf{z}_2'}} \|\mathbf{z}-\mathbf{c}_B\|.
\end{align}
Given that the distance between the line segment and point is determined by the shortest path, the distance between $\overline{\mathbf{z}_1\mathbf{z}_2}$ and $\mathbf{c}_B$ is equal to the distance between the chord $\overline{\mathbf{z}_1'\mathbf{z}_2'}$ and $\mathbf{c}_B$, i.e., $d(\overline{\mathbf{z}_1\mathbf{z}_2},\mathbf{c}_B) =d(\overline{\mathbf{z}_1'\mathbf{z}_2'},\mathbf{c}_B)$.

Let $l'$ be the length of a chord $\overline{\mathbf{z}_1'\mathbf{z}_2'}$; then, $d(\overline{\mathbf{z}_1'\mathbf{z}_2'},\mathbf{c}_B)$ is determined by
\begin{align}
    d(\overline{\mathbf{z}_1'\mathbf{z}_2'},\mathbf{c}_B) = \sqrt{\mathrm{R}_B^2 - \left(\frac{l'}{2}\right)^2}.
\end{align}
Because $d(\overline{\mathbf{z}_1'\mathbf{z}_2'},\mathbf{c}_B)$ is a decreasing function of $l'$ for $0 < l' \leq L$, the following relationship can be established. 
\begin{align}
    \min_{0< l' \leq L} d(\overline{\mathbf{z}_1'\mathbf{z}_2'},\mathbf{c}_B) = \sqrt{\mathrm{R}_B^2 \!-\!\left(\!\frac{L}{2}\!\right)^{\!2}} \leq \sqrt{\mathrm{R}_B^2 \!-\!\left(\!\frac{l'}{2}\!\right)^{\!2}}. 
\end{align}

Finally, we can conclude that $d(\overline{\mathbf{z}_1\mathbf{z}_2},\mathbf{c}_B)$ is greater than or equal to $\sqrt{\mathrm{R}_B^2 -\left(\frac{L}{2}\right)^2}$. \qed


\end{document}